\documentclass[11pt,reqno]{article}
\usepackage{jmlr2e}
\usepackage{fullpage,times,graphicx,amssymb,amsmath,amsfonts,bbm,psfrag,xcolor}

\usepackage{graphicx}
\usepackage{amssymb}
\usepackage{amsmath}
\usepackage{amsfonts}
\usepackage{bbm}
\usepackage{psfrag}
\usepackage{xcolor}
\usepackage{fullpage}
\usepackage{epstopdf}
\usepackage{graphicx}
\usepackage{times}
\usepackage{jmlr2e}
\usepackage{commath}
\usepackage{thmtools}
\usepackage{thm-restate}
\usepackage{times}
\usepackage{graphicx}
\usepackage{natbib}
\usepackage{dirtytalk}
\usepackage{empheq}
\definecolor{ddarkbrown}{rgb}{0.5,0.2,0.05} \definecolor{bbluegray}{rgb}{0.05,0,0.5}

\usepackage[colorlinks,citecolor=bbluegray,linkcolor=ddarkbrown,urlcolor=blue,breaklinks]{hyperref}

\usepackage{subfig,float} 

\usepackage{mathtools}

\usepackage{algorithm,algcompatible,algpseudocode}
\usepackage{multicol,lipsum}
\algnewcommand{\Inputs}[1]{%
	\State \textbf{Inputs: \:}{#1}
}

\algnewcommand{\Output}[1]{%
	\State \textbf{Output: \:}{#1}
}
\algnewcommand{\Initialize}[1]{%
	\State \textbf{Initialize: \:}{#1}
}

\algnewcommand{\IIf}[1]{\State\algorithmicif\ #1\ \algorithmicthen}
\algnewcommand{\EndIIf}{\unskip\ \algorithmicend\ \algorithmicif}

\let \oldsection \section
\renewcommand{\section}{\vspace{3ex plus 1ex}\oldsection}

\newcommand{\BEAS}{\begin{eqnarray*}}
	\newcommand{\EEAS}{\end{eqnarray*}}
\newcommand{\BEA}{\begin{eqnarray}}
\newcommand{\EEA}{\end{eqnarray}}
\newcommand{\mb}{\mathbb}
\newcommand{\BEQ}{\begin{equation}}
\newcommand{\EEQ}{\end{equation}}
\newcommand{\BIT}{\begin{itemize}}
	\newcommand{\EIT}{\end{itemize}}
\newcommand{\BNUM}{\begin{enumerate}}
	\newcommand{\ENUM}{\end{enumerate}}

    \newcommand{\Pk}{\mathcal{P}_{\mathcal{K}}}
    
    \newcommand{\Pc}{\mathcal{P}_{\mathcal{C}}}

    \newcommand{\Gr}{G^\star}

	\newcommand{\E}{\mathbb{E}}	

\newcommand{\vi}{v \in \mathcal{C}\cap\mathcal{B}^d}
\newcommand{\Agi}{A_{g_i}}
\newcommand{\Tgi}{T_{g_i}}

	\newcommand{\ve}{\varepsilon}
\newcommand{\BA}{\begin{array}}
	\newcommand{\EA}{\end{array}}

\newcounter{dummy} \numberwithin{dummy}{section}
\newtheorem{mythm}[dummy]{Theorem}

\newtheorem{mylem}[dummy]{Lemma}

\numberwithin{mythm}{section}
\numberwithin{mydef}{section}
\numberwithin{myprop}{section}
\numberwithin{mylem}{section}
\numberwithin{mycor}{section}

\title{Group Symmetry Enables Faster Optimization in Inverse Problems}

\begin{document}	

	\author{\name Junqi Tang  \email j.tang.2@bham.ac.uk\\
		\addr School of Mathematics,\\ University of Birmingham \\
        \\
        \name Guixian Xu  \email gxx422@student.bham.ac.uk\\
		\addr School of Mathematics,\\ University of Birmingham \\
        }

	\editor{}

	\maketitle


\begin{abstract}

We prove for the first time that, if a linear inverse problem exhibits a group symmetry structure, gradient-based optimizers can be designed to exploit this structure for faster convergence rates. This theoretical finding demonstrates the existence of a special class of structure-adaptive optimization algorithms which are tailored for symmetry-structured inverse problems such as CT/MRI/PET, compressed sensing, and image processing applications such as inpainting/deconvolution, etc.

\end{abstract}

\section{Introduction}
In this work we consider linear inverse problems of the form:
\begin{equation}
    b = A x^\dagger + w, A \in \mathbb{R}^{m \times d}
\end{equation}
where $x^\dagger \in \mathcal{X} \subseteq \mb{R}^d$ being the ground-truth signal/image to be estimated, $A$ being the observation forward operator, $y$ being the observation, and $w$ being the measurement noise. Given a constraint set $\mathcal{K}\subseteq \mathcal{X}$, one can obtain a classical least-square estimator of the form:
\begin{equation}
    x^\star \in \arg\min_{x \in \mathcal{K}} f(x):=\frac{1}{2}\|Ax - b\|_2^2.
\end{equation}
The standard first-order solver for this optimization problem would be the proximal (projected) gradient descent \citep{combettes2011proximal} method:
\begin{equation}
    x_{k+1} = \mathcal{P}_{\mathcal{K}}[x_k - \eta \nabla f(x_k)],
    \tag{PGD}
\end{equation}
where $\Pk$ is the proximal operator on the indicator function $\iota_\mathcal{K}$ of the set $\mathcal{K}$ (orthongal projection on the set $\mathcal{K}$):
\begin{equation}
    \Pk(x) = \arg\min_{y \in \mb{R}^d} \|x - y\|_2^2 + \iota_{\mathcal{K}}(y) = \arg\min_{y \in \mathcal{K}} \|x - y\|_2^2
\end{equation}
To effectively analyze the convergence behavior of projected gradient descent (PGD) for the difficult scenarios where strong-convexity of the objective is not available (which is typically the case for inverse problems), or the constraint is non-convex, \cite{oymak2017sharp} have developed an accurate theoretical framework which establishes sharp bounds for the linear convergence rate PGD algorithm under restricted strong-convexity \citep{agarwal2012fast}:
\begin{equation}
    \|Av\|_2^2 \geq \mu_C\|v\|_2^2,\ \ \forall v \in \mathcal{C}_{\mathcal{K}-x^\dagger}
\end{equation}
where $\mathcal{C}_{\mathcal{K}-x^\dagger}$ being the descent cone at $x^\dagger$ covering the shifted set $\mathcal{K}-x^\dagger$:
\begin{equation}
    \mathcal{C}_{\mathcal{K}-x^\dagger} := \{v \in \mb{R}^d | v=a(x - x^\dagger), \forall a \geq 0, x \in \mathcal{K}\}
\end{equation}
Let $L = \|A^TA\|$ (largest eigen-value of the Hessian $A^TA$), a linear convergence rate can be established for PGD as demonstrated by \cite{oymak2017sharp} despite the lack of standard strong-convexity:
\begin{equation}
    \|x_{k+1} - x^\dagger\|_2 \leq \alpha^{k/2} \|x_0 - x^\dagger\|_2 + O(\|w\|_2)
\end{equation}
where $\alpha = \kappa_c(1 - \frac{\mu_C}{L})^{\frac{1}{2}}$. This result demonstrates that the intrinsic low-dimensional structure of the solution can be exploited to achieve faster convergence of the optimizer if the forward operator is well suited for the signal set \citep{pmlr-v70-tang17a}.

More recently, researchers have eventually determined that the spectral structure of the measurement forward operator can be utilized in optimization algorithms for faster convergence, this is, in fact, the precise reason for the success of stochastic gradient-based methods in solving inverse problems such as CT / PET \citep{tang2020practicality,ehrhardt2024guide}.

In this work, we explore from a theoretical point of view whether a gradient-based optimizer can utilize a more delicate structure, namely the group symmetry structure \citep{tachella2023sensing}, for faster convergence. Even in the case where restricted strong convexity is unavailable ($\mu_C = 0$), the forward operator in practice can often be viewed as a subsampling of a complete measurement $A_{full}$ -- for example, the sparse-view or limited-angle CT operator can be viewed as a subsampling of a full-angle CT scan. If this extra information not utilized, the PGD cannot have a linear convergence but only a slow rate on the objective function:
\begin{equation}
    f(x_K) - f(x^\star) \leq O(1/K),
\end{equation}
or $O(1/K^2)$ at best for FISTA / Nesterov's acceleration \citep{2009_Beck_Fast,nesterov2007gradient}, both without any guarantee for the convergence rate towards the ground-truth. Here we show that this rate can still be improved to a linear rate similar to the one proven by \cite{oymak2017sharp} but with a relaxed condition on restricted strong convexity utilizing a group-symmetry structure.

\section{Hunting in the shadow: exploiting the group symmetry structure for acceleration}

In practice, many measurement operators can be viewed as a subsampling of a full operator. Given a cyclic group $(G, \circ)$ which is generated by an element $g \in G$:
\begin{equation}
    G = \langle g \rangle,
\end{equation}
with group actions on the signal set $\mathcal{X}$:
\begin{equation}
    \Tgi : G \times \mathcal{X} \rightarrow \mathcal{X}
\end{equation}
Typical examples include rotations and translations. For rotations (typically for medical image tasks such as CT/MRI), suppose we index the pixels of an image in terms of radians $r$ and angles $\theta$ \citep{tachella2023sensing}:
\begin{equation}
    \Tgi x = x(r, \theta - g_i)
\end{equation}
Then we can derive new measurement operators in addition to $A$, but without actually having the corresponding measurement data for them. Let $g_1 := g$, $g_2 := g\circ g = g^2$, $g_3 := g^3, ...$, we have:
\[
A_{\mathrm{full}} :=\begin{bmatrix}
A \\
A_{g_1} \\
\vdots \\
A_{g_{|G|-1}}
\end{bmatrix}=\begin{bmatrix}
A \\
AT_{g_1} \\
\vdots \\
AT_{g_{|G|-1}}
\end{bmatrix}
\]
Let's use X-ray CT for a convenient illustrative example here. Suppose $A$ consists of X-ray measurements from a subset of angles in $[0, 360]$ degrees, while $g$ being the rotation of 1 degree, then $A_{\mathrm{full}}$ will cover all 360 degrees with the group size being $|G|=360$ (considering the subgroup containing only the integer degrees), leading to a much better conditioning and satisfying the restricted strong-convexity condition much easier than $A$ alone. To utilize $A_{\mathrm{full}}$ we can modify the PGD with group actions as such:

\begin{equation}
    x_{k+1} = \mathcal{P}_{\mathcal{K}}[x_k - \eta \Tgi^{-1} \nabla f(\Tgi x_k)],\ \ \mathrm{Sample\ } g_i \sim G, \mathrm{\ s.t.\ a\  distribution\ } P
    \tag{Group-PGD}
\end{equation}
which is essentially: $x_{k+1}= \mathcal{P}_{\mathcal{K}}[x_k - \eta \Agi^T (\Agi x_k - b)]$ for linear regression case. We name this algorithm Group-PGD. In each of iteration of the Group-PGD, a group action $\Tgi$ is randomly selected according to some distribution $P$ to perturb first the current iterate $x_k$, then the gradient $\nabla f$ is evaluated at the perturbed point $\Tgi x_k$, followed by $\Tgi^{-1}$. 

Intuitively, we should focus the sampling distribution on those $\Tgi$'s that do not deviate too much from the image, that is, $\|\Tgi x^\dagger - x^\dagger\|_2 \approx 0$, which can also be observed via our theory. For instance, in terms of rotations in sparse-view CT reconstruction tasks, we found that we should choose $g_i$ to be mild rotations of $\pm1, \pm2$ degrees, which would suffice for us to observe significant acceleration over standard PGD. We denote this subset of $G$ as $G_\star$ and make the sampling distribution $P$ supported only on $G_\star$. This symmetric subset $G_\star \subseteq G$ does not need to form a subgroup but should include the identity element and the inverses of each element (just satisfying the identity axiom and the inverse axiom, not the closure axiom). For this case, we define an complemented forward operator that shall be better conditioned compared to $A$ alone:
\begin{equation}
A_{G_\star} :=\frac{1}{|G_\star|}\begin{bmatrix}
A \\
A_{g_1} \\
\vdots \\
A_{g_{|G_\star|-1}}
\end{bmatrix}=\frac{1}{|G_\star|}\begin{bmatrix}
A \\
AT_{g_1} \\
\vdots \\
AT_{g_{|G_\star|-1}}
\end{bmatrix}
\label{agi}
\end{equation}
where $g_1 = g$, $g_2 = g^{-1}$, $g_3 = g^2$, $g_4 = g^{-2}, .....$ With $A_{G_\star}$ we can build convergence theorem for Group-PGD under the assumption of restricted strong-convexity with $A_{G_\star}$:
\begin{equation}
    \|A_{G_\star} v\|_2 \geq \mu_{G_\star}\|v\|_2,\ \ \forall v \in \mathcal{C}_{\mathcal{K}-x^\dagger}
\end{equation}
instead of plain restricted strong convexity on $A$ (for extreme undetermined cases $\mu_C \approx 0$):
\begin{equation}
    \|A v\|_2 \geq \mu_{C}\|v\|_2,\ \ \forall v \in \mathcal{C}_{\mathcal{K}-x^\dagger}
\end{equation}
with $\mu_{G_\star} \gg \mu_C$ for ill-posed problems. For our theoretical analysis, we will utilize the following results from the literature:

\begin{mylem}
    (Projection identities, \citep{oymak2017sharp}) Given a set $\mathcal{K}$, $\mathcal{B}^d$ the unit ball in $\mb{R}^d$ and an orthogonal projection operator $\Pk$, $x^\dagger \in \mathcal{K}$ and 
    \begin{equation}
    \mathcal{C} := \{v \in \mb{R}^d | v=a(x - x^\dagger), \forall a \geq 0, x \in \mathcal{K}\}
\end{equation}
    we have:
    \begin{equation}
        \|\Pc(x)\|_2 = \sup_{v \in \mathcal{C} \cap \mathcal{B}^d} \langle v, x\rangle,
    \end{equation}
    and
    \begin{equation}
        \Pk(x+v)-x = \mathcal{P}_{\mathcal{K}-x}(v),
    \end{equation}
    and
    \begin{equation}
        \|\Pk(x)\|_2 \leq \kappa_c\|\Pc(x)\|_2
    \end{equation}
    where $\kappa_c = 1$ if $\mathcal{K}$ is convex, $\kappa_c = 2$ if $\mathcal{K}$ is nonconvex.
\end{mylem}


\section{Main result}

\begin{mythm}
    Given measurements $b = Ax^\dagger + w$. Let $\eta = \frac{1}{L}$ where $L= \|A^TA\|$, and $G_\star$ is a symmetric subset of $G$ including identity, while $A_{G_\star}$ defined in \eqref{agi}, suppose:
    \begin{equation}
    \|A_{G_\star} v\|_2 \geq \mu_{G_\star}\|v\|_2,\ \ \forall v \in \mathcal{C}_{\mathcal{K}-x^\dagger}
\end{equation}
    then the following bound holds for Group-PGD (at $K$-th iteration) if for each iteration $g_i$ is sampled from $G_\star$ uniformly at random:
\begin{equation}
     \E\|x_{K+1} - x^\dagger\|_2 \leq \alpha_{\Gr}^K \|x_0 - x^\dagger\|_2 + \frac{\kappa_c(1- \alpha_{\Gr}^K)}{L(1-\alpha_{\Gr})}(\ve_{\Gr} + \ve_w\|w\|_2)
\end{equation}
where $\kappa_c = 1$ if $\mathcal{K}$ is convex, $\kappa_c = 2$ if $\mathcal{K}$ is nonconvex, and:
\begin{equation}
    \alpha_{\Gr} := \kappa_c  (1- \frac{\mu_{\Gr}}{L})^\frac{1}{2}
\end{equation}
\begin{equation}
    \ve_{\Gr} := \sup_{\vi,g_i\in \Gr} v^T\Agi^TA(x^\dagger - \Tgi x^\dagger),
\end{equation}
\begin{equation}
    \ve_w := \sup_{\vi,g_i\in \Gr} v^T\Agi^T\frac{w}{\|w\|_2}
\end{equation}
\end{mythm}

\begin{remark}
    Compared to the PGD convergence bound proven by \cite{oymak2017sharp}, the above bound for Group-PGD enjoys a much faster linear rate dependent on $\mu_{G_\star}$, at a price of the error term $\ve_{\Gr}$. To have a closer look:
    \begin{equation}
        \ve_{\Gr} := \sup_{\vi,g_i\in \Gr} v^T\Agi^TA(x^\dagger - \Tgi x^\dagger) \leq \max_{g_i\in \Gr} \|\Agi^TA\|\|x^\dagger - \Tgi x^\dagger\|_2
    \end{equation}
    To make the term small, one must carefully choose the symmetric subset $G_\star$ such that $\|\Tgi x^\dagger - x^\dagger\|_2$ is small. For instance, rotations of $\pm1, \pm2$ degrees for CT reconstruction.
\end{remark}

\begin{remark}
    In practice, considering this theoretical result, it should be advised that a multistage scheme which starts with a large symmetric set and then recursively shrinks its size should be applied for best performance. We leave the study of this speed-accuracy trade-off for future work.
\end{remark}

\subsection{The proof of Theorem 3.1}
We present the proof of our main theorem here:

\begin{proof}

Recall the iterations of Group-PGD:
\begin{equation}
    x_{k+1} = \mathcal{P}_{\mathcal{K}}[x_k - \eta T_{g_i}^{-1} \nabla f(T_{g_i}x_k)] = \mathcal{P}_{\mathcal{K}}[x_k - \eta \Agi^T (\Agi x_k - b)]
\end{equation}
we can bound the estimation error $\|x_{k+1} - x^\dagger\|_2$ using Lemma 2.1, and several manipulations:
\begin{eqnarray*}
    \|x_{k+1} -x^\dagger\|_2 &=& \|\Pk[x_k - \eta \Agi^T (\Agi x_k - b)]  - x^\dagger\|_2 \\
    &=&\|\mathcal{P}_{\mathcal{K}-x^\dagger}[x_k - \eta \Agi^T (\Agi x_k - b) - x^\dagger]\|_2\\
     &\leq&\kappa_c\|\Pc[x_k - x^\dagger - \eta \Agi^T (\Agi x_k - b)]\|_2\\
          &\leq&\kappa_c\|\Pc[x_k - x^\dagger - \eta \Agi^T (\Agi x_k - Ax^\dagger - w)]\|_2\\
        &\leq&\kappa_c\|\Pc[x_k - x^\dagger - \eta \Agi^T \Agi (x_k-x^\dagger) + \eta \Agi^TA(x^\dagger - \Tgi x^\dagger) + \eta\Agi^Tw)]\|_2\\
        &\leq&\kappa_c\sup_{\vi}v^T[x_k - x^\dagger - \eta \Agi^T \Agi (x_k-x^\dagger) + \eta \Agi^TA(x^\dagger - \Tgi x^\dagger) + \eta\Agi^Tw)]\\
    &\leq&\kappa_c\sup_{\vi}v^T[I- \eta \Agi^T \Agi] (x_k-x^\dagger) + \eta\kappa_c\sup_{\vi} v^T[\Agi^TA(x^\dagger - \Tgi x^\dagger) + \eta\Agi^Tw)]\\
    &\leq&\kappa_c\|(I- \eta \Agi^T \Agi)(x_k-x^\dagger)\|_2 
\\ &&+ \eta\kappa_c\sup_{\vi} v^T\Agi^TA(x^\dagger - \Tgi x^\dagger) + \eta\kappa_c\sup_{\vi} v^T\Agi^Tw\\
\end{eqnarray*}
Then taking expectation w.r.t. sampling distribution $P$, and then apply Jensen's inequality:
\begin{eqnarray*}
    \E\|x_{k+1} - x^\dagger\|_2 &\leq& \kappa_c \E\|(I- \eta \Agi^T \Agi)(x_k-x^\dagger)\|_2 \\ &&+ \eta\kappa_c\sup_{\vi,g_i\in \Gr} v^T\Agi^TA(x^\dagger - \Tgi x^\dagger) + \eta\kappa_c\sup_{\vi,g_i\in \Gr} v^T\Agi^Tw\\
    &\leq& \kappa_c \sqrt{\E\|(I- \eta \Agi^T \Agi)(x_k-x^\dagger)\|_2^2} \\ &&+ \eta\kappa_c\sup_{\vi,g_i\in \Gr} v^T\Agi^TA(x^\dagger - \Tgi x^\dagger) + \eta\kappa_c\sup_{\vi,g_i\in \Gr} v^T\Agi^Tw\\
        &\leq& \kappa_c \{\E[\|(x_k-x^\dagger)\|_2^2 - 2\eta\|\Agi(x_k - x^\dagger)\|_2^2 + \eta^2\|\Agi^T\Agi(x_k - x^\dagger)\|_2^2)\}^{0.5} \\ &&+ \eta\kappa_c\sup_{\vi,g_i\in \Gr} v^T\Agi^TA(x^\dagger - \Tgi x^\dagger) + \eta\kappa_c\sup_{\vi,g_i\in \Gr} v^T\Agi^Tw\\ 
\end{eqnarray*}
Since $\Tgi$ are all unitary transforms, and $\eta < \frac{2}{L}$, we will have:
\begin{eqnarray*}
    \E\|x_{k+1} - x^\dagger\|_2 
        &\leq& \kappa_c \{\E[\|(x_k-x^\dagger)\|_2^2 - 2\eta\|\Agi(x_k - x^\dagger)\|_2^2 + \eta^2L\|\Agi(x_k - x^\dagger)\|_2^2)\}^{0.5} \\ &&+ \eta\kappa_c\sup_{\vi,g_i\in \Gr} v^T\Agi^TA(x^\dagger - \Tgi x^\dagger) + \eta\kappa_c\sup_{\vi,g_i\in \Gr} v^T\Agi^Tw\\
        &\leq& \kappa_c \{\E[\|(x_k-x^\dagger)\|_2^2 - \eta(2 -\eta L)\|\Agi(x_k - x^\dagger)\|_2^2\}^{0.5} \\ &&+ \eta\kappa_c\sup_{\vi,g_i\in \Gr} v^T\Agi^TA(x^\dagger - \Tgi x^\dagger) + \eta\kappa_c\sup_{\vi,g_i\in \Gr} v^T\Agi^Tw\\
        &\leq& \kappa_c \{\|(x_k-x^\dagger)\|_2^2 - \eta(2 -\eta L)\mu_G\|(x_k - x^\dagger)\|_2^2\}^{0.5} \\ &&+ \eta\kappa_c\sup_{\vi,g_i\in \Gr} v^T\Agi^TA(x^\dagger - \Tgi x^\dagger) + \eta\kappa_c\sup_{\vi,g_i\in \Gr} v^T\Agi^Tw\\
        &\leq& \kappa_c \{ [1- \eta(2 -\eta L)\mu_{\Gr}]\|(x_k - x^\dagger)\|_2^2\}^{0.5} \\ &&+ \eta\kappa_c\sup_{\vi,g_i\in \Gr} v^T\Agi^TA(x^\dagger - \Tgi x^\dagger) + \eta\kappa_c\sup_{\vi,g_i\in \Gr} v^T\Agi^Tw\\
\end{eqnarray*}
Taking $\eta = \frac{1}{L}$, we have:
\begin{eqnarray*}
    \E\|x_{k+1} - x^\dagger\|_2 &\leq& \kappa_c  (1- \frac{\mu_{\Gr}}{L})^\frac{1}{2}\|x_k - x^\dagger\|_2 \\ &&+ \eta\kappa_c\sup_{\vi,g_i\in \Gr} v^T\Agi^TA(x^\dagger - \Tgi x^\dagger) + \eta\kappa_c\sup_{\vi,g_i\in \Gr} v^T\Agi^Tw\\
\end{eqnarray*}
Let $\alpha_{\Gr} = \kappa_c  (1- \frac{\mu_{\Gr}}{L})^\frac{1}{2}$, we can tower-up the recursion to $x_0$:
\begin{eqnarray*}
    \E\|x_{k+1} - x^\dagger\|_2 &\leq& [\kappa_c  (1- \frac{\mu_{\Gr}}{L})]^\frac{k}{2}\|x_0 - x^\dagger\|_2 \\ &&+ \frac{\kappa_c(1- \alpha_{\Gr}^k)}{L(1-\alpha_{\Gr})}\sup_{\vi,g_i\in \Gr} v^T\Agi^TA(x^\dagger - \Tgi x^\dagger) \\&&+ \frac{\kappa_c(1- \alpha_{\Gr}^k)\|w\|_2}{L(1-\alpha_{\Gr})}\sup_{\vi,g_i\in \Gr} v^T\Agi^T\frac{w}{\|w\|_2}\\
\end{eqnarray*}
\end{proof}
\clearpage
\begin{figure}
    \centering
    \includegraphics[width=0.95\linewidth]{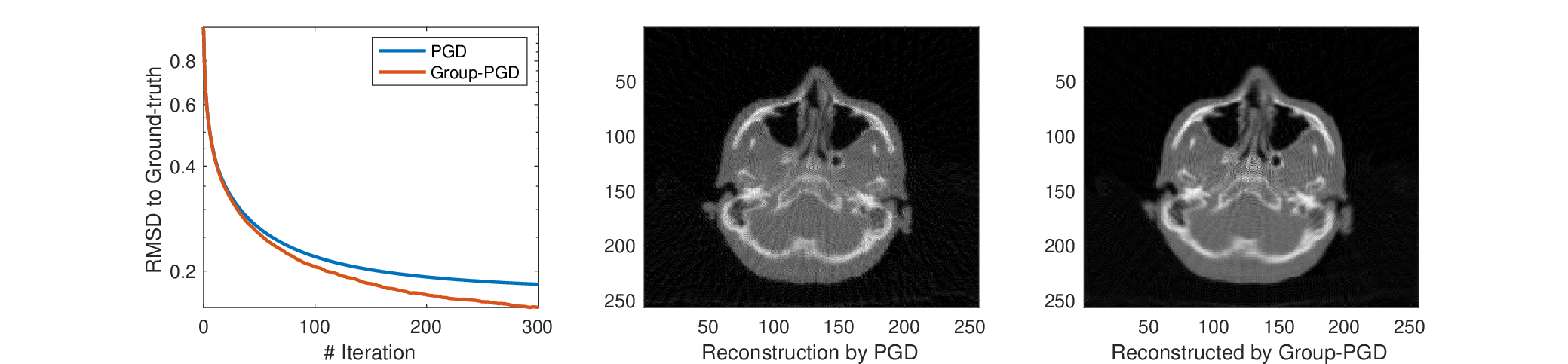}
    \caption{Sparse-view fan-beam CT example, comparing PGD with Group-PGD}
    \label{fig1}
\end{figure}

\begin{figure}
    \centering
    \includegraphics[width=0.95\linewidth]{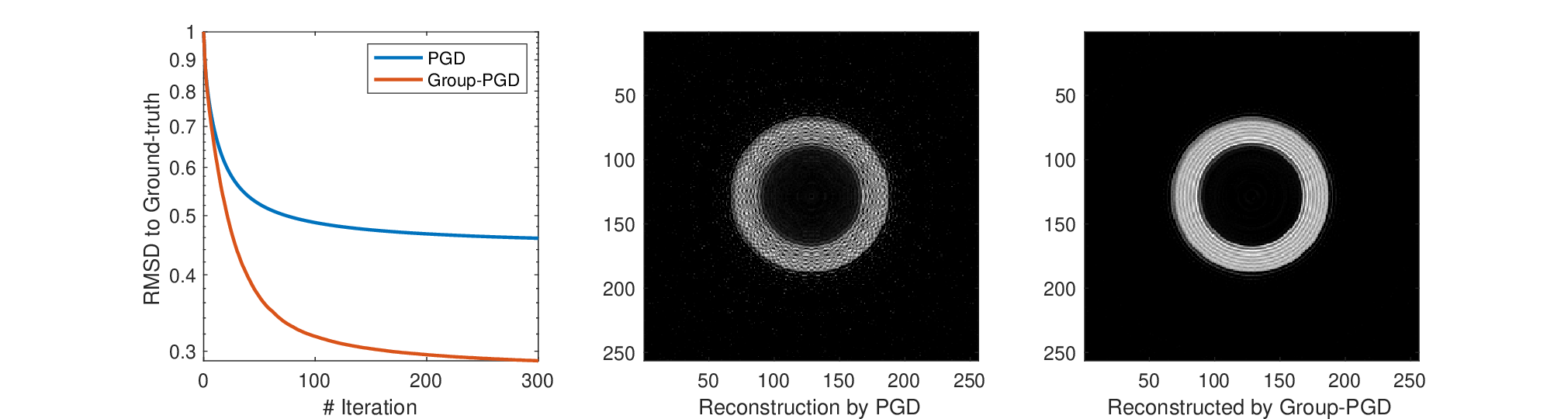}
    \caption{Extreme sparse-view fan-beam CT example, comparing PGD with Group-PGD}
    \label{fig2}
\end{figure}

\section{Numerical verification}

In this section, we present some numerical results on Group-PGD compared to standard PGD on sparse-view fan-beam CT. For a proof of concept, we consider the simplest case where only a box-constraint is enforced ($b$ here is measurement data corrupted by Poisson noise):
\begin{equation}
    x^\star \in \arg\min_{x \in [0,1]^d} f(x):=\frac{1}{2}\|Ax - b\|_2^2,
\end{equation}
hence clearly the standard restricted strong-convexity parameter $\mu_C = 0$, where PGD is expected to have a slow convergence rate. For Group-PGD, we can utilize the group transforms to make $\mu_{G_\star} > 0$ and activate fast convergence. Here we consider two examples, the first on a real CT image and the second one on a generated simple ring image. 

For the first example in Figure \ref{fig1}, we choose $|G_\star| = 5$ where $G_\star =\{\mathrm{Id}, g, g^2, g^{-1}, g^{-2}\}$ (note that this is not a true subgroup of $G$ as it does not satisfy the closure axiom) and $g$ is a one-degree rotation. Since this is a complicated image to reconstruct we have to choose small degrees of rotation to ensure $\Tgi x^\dagger \approx x^\dagger$ hence $\ve_{G_\star}$ is small. For this example, we have $A \in \mathbb{R}^{22344 \times 65536}$, which the number of measurements is around $34\%$ of the number of pixels (unknowns). The operator $A_{G_\star}$ will include measurement physics from all 360 degrees and has a number of measurements $170\%$ of the number of pixels. We can observe a significant acceleration of Group-PGD over standard PGD in terms of convergence rate in root mean square distance (RMSD) towards ground truth.

For the second example in Fig. \ref{fig2} we tested on a ring image which satisfies $\Tgi x^\dagger = x^\dagger$ for arbitrary rotations and hence $\ve_{G_\star} = 0$. For this example, we have an extreme sparse-view case which $A \in \mathbb{R}^{4104 \times 65536}$, in which the number of measurement is around $6\%$ of the number of pixels. For this case since the image is simple, we can be more greedy on $G_\star$ (we choose $|G_\star| = 55$ here). We can observe a huge acceleration in this case as a proof of concept.

\section{Conclusion}

In this work, we present a fundamental algorithmic design of symmetric-adaptive first-order methods for linear inverse problems and a convergence proof showing that a fast linear convergence rate can be achieved when standard restricted strong-convexity is not available, if the algorithm utilize the group-symmetry structure of the inverse problem. We believe that this work justifies the motivation of a new research direction for optimization algorithms in inverse problems, extending previous studies which consider utilizing only the intrinsic low-dimensionality structures (such as sparsity/low-rankness, etc.) in solutions and measurement operators for acceleration using randomized sketching or stochastic gradient-based optimization \citep{driggs2021stochastic,ehrhardt2024guide}. Meanwhile in the line of work on plug-and-play algorithms \citep{terris2024equivariant,tang2024practical}, group transforms have been recently applied to improve the stability and performance of deep denoisers applied in PnP as off-the-shelf deep denoisers often lead to instability (Lipschitz constant greater than 1). It will also be an interesting direction to discover the deeper reason for the success of this equivariant denoiser and the implications of it in terms of optimization.

\bibliography{main.bib}

\end{document}